\documentclass{amsart}
\usepackage{latexsym}
\usepackage{amsfonts}
\usepackage{amssymb,mathrsfs,calrsfs}
\usepackage{amscd}
\usepackage{epsfig}

\newtheorem{theorem}{Theorem}

\newtheorem{kor}{Corollary}
\newtheorem{prop}{Proposition}

\newlength{\figboxwidth}
\newcommand{\lat}{\Z \oplus \Z i}

\newcommand{\abx}{(X,\omega)}
\newcommand{\aby}{(Y,\tau)}
\newcommand{\abas}{(S,\alpha)}

\newcommand{\slr}{ {\rm SL}_2(\R) }

\newcommand{\slz}{ {\rm SL}_2(\Z) }
\newcommand{\pslz}{ {\rm PSL}_2(\Z) }

\newcommand{\slv}[1]{ {\rm SL}(#1) }

\newcommand{\R}{\mathbb{R}}
\newcommand{\C}{\mathbb{C}}
\newcommand{\Q}{\mathbb{Q}}
\newcommand{\Z}{\mathbb{Z}}

\newcommand{\proj}{\mathbb{CP}}
\newcommand{\T}{\mathbb{T}}
\renewcommand{\H}{\mathbb{H}}

\newcommand{\ct}{\mathfrak{t}}

\newcommand{\sC}{\mathscr{C}}
\newcommand{\sE}{\mathscr{E}}
\newcommand{\sF}{\mathscr{F}}

\newcommand{\sL}{\mathscr{L}}
\newcommand{\sM}{\mathscr{M}}
\newcommand{\sO}{\mathscr{O}}

\newcommand{\CF}{\mathcal{F}}

\newcommand{\CH}{\mathcal{H}}

\newcommand{\CS}{\mathcal{S}}

\newcommand{\CU}{\mathcal{U}}

\newcommand{\ca}{\mathfrak{a}}
\newcommand{\cb}{\mathfrak{b}}

\DeclareMathOperator{\area}{area}

\DeclareMathOperator{\aut}{Aut}

\DeclareMathOperator{\D}{D}

\DeclareMathOperator{\id}{id}

\DeclareMathOperator{\per}{Per}
\DeclareMathOperator{\hol}{hol}

\DeclareMathOperator{\pr}{pr}
\DeclareMathOperator{\im}{Im}

\title{Spaces of elliptic differentials    
}
\author{Martin Schmoll}
\date{\today}

\begin{document}

\begin{abstract}
We study modular fibers of elliptic differentials, i.e.  
investigate spaces of coverings $\aby \rightarrow (\C/\lat,dz)$.  
For genus $2$ torus covers with fixed degree 
we show, that the modular fibers $\sF_d(1,1)$ 
are connected torus covers with Veech group $\slz$. 
Using results of Eskin, Masur and Schmoll \cite{ems} 
we calculate $\chi(\sF_d(1,1))$ and  
the parity of the spin structure of the 
quadratic differential $(\sF_d(1,1)/(-\id),q_d)$.  
We state and apply formul{\ae} for the asymptotic 
quadratic growth rates of various types of geodesic 
segments on $\aby \in \sF_d(1,1)$. The quadratic growth 
rates are expressed in terms of the $\slz$ 
orbit closure of $\aby$ in $\sF_d(1,1)$ and the flat geometry of 
$\sF_d(1,1)$.\\
These are extended notes from a talk the author gave 
during the {\em Activity on Algebraic and Topological Dynamics} at 
the Max-Planck-Institute for Mathematics, Bonn summer 2004.   
\end{abstract}
\address{
 The Pennsylvania State University,  
 7E Thomas, 
 University Park, PA, 16802 }
\email{schmoll@math.psu.edu}
\maketitle
\tableofcontents


\section{Introduction}

{\bf Motivation and Background.} If we want to find the length 
distribution of isotopy classes 
of closed geodesics on the flat torus $\T^2:=\C/\lat \cong \R^2/\Z^2$, 
the answer is easily obtained by counting integer lattice points 
in the plane: 
\[N(\T^2,T):=|\{(x,y) \in \Z^2: \ \gcd(x,y)=1, \sqrt{x^2+y^2}<T \}|
\sim \frac{\pi}{\zeta(2)}T^2=\frac{6}{\pi}T^2.\] 
The factor $\frac{1}{\zeta(2)}$ arises, 
if one counts {\em primitive} geodesics (see \cite{em}) including 
their direction, counting of geodesics ignoring direction requires 
the weight $\frac{1}{2\zeta(2)}$. 
Primitive geodesics (with direction) 
are represented by integer lattice points in $\R^2$ 
which are {\em visible from the origin}.
We like to ask the same question for a (branched) covering 
$\pi: X \rightarrow \C/\Lambda$, $\Lambda \subset \C$ 
a lattice.   
The necessary (complex) geometric structure on $X$ is a 
{\em holomorphic differential} $\omega$ obtained by pulling back 
the differential $dz$ on $\C/\Lambda$. We call the pair 
$(X,\omega =\pi^{\ast}dz)$ {\em elliptic differential}. 
Using $\omega$ one 
identifies $X$ locally with regions in the complex plane 
by coordinates of the shape 
\[ z_{p_0}(p)=\int^{p}_{p_0}\omega\] 
away from the zero-set $Z(\omega)$ of $\omega$. 
With respect to these charts 
coordinate changes are translations, in particular the 
Euclidean metric pulls back to $(X,\omega)$ and defines a 
global Euclidean metric on $X-Z(\omega)$. 
\medskip\\
We are mainly interested in the following geodesic segments 
on $\abx$: 
\begin{itemize}
\item (isotopy classes of) closed geodesics $Cyl(\omega)$ and 
\item saddle connections $SC(\omega)$, these are geodesic 
segments starting and ending at zeros of $\omega$, 
without hitting a zero in between.      
\end{itemize}
We like to study the asymptotic quadratic 
growth rate of these geodesic segments on $\abx$ 
with respect to the Euclidean metric defined by $\omega$, 
i.e. we look at the number 
\begin{equation}
N_{cyl}(\omega,T):=\left| \left\{\gamma \in Cyl(\omega): 
\int_{\gamma}|\omega|<T \right\}\right|.
\end{equation}
A fundamental result of Masur 
(for a new version see \cite{em}) says
\begin{theorem}\cite{m3,m4}
For any translation surface $\abx$ there are constants, 
such that for $T>>0$
\[0 < c_{1}T^2 < N_{cyl}(\omega,T)<c_{2}T^2.\]
The same is true for the set of saddle connections $SC(\omega)$ 
on $\abx$ (eventually with different constants $c_i$).
\end{theorem} 
Surprisingly in various cases \cite{v2,ems,emz,mcm3,emm}, 
there is an asymptotic quadratic formula 
\[N_{cyl}(\omega,T)\sim \frac{\pi}{\zeta(2)}c_{cyl}(\omega)T^2. \]
Moreover: in all cases known so far the constant 
can be computed \cite{v3,ems}, or at least expressed in terms of 
geometrical data of the moduli space of 
Abelian differentials where $\abx$ belongs too \cite{emz}.
For general differentials $\abx$ it is not known that the various 
asymptotic constants $c_{\ast}(\omega)$ exist. In 
the case of elliptic differentials however, it is well known 
(see \cite{ems,emm,s1}) that all asymptotic constants, including   
$c_{cyl}(\omega)$ and $c_{SC}(\omega)$  are well defined. 

To introduce the main objects we  
give an alternative, more geometric description of translation 
surfaces.\medskip\\
{\bf Translation surfaces by gluing polygons.} 
Take a finite set of polygons $P_1,...,P_n$ in the complex plane $\C$ 
with boundary components $\partial P_i$ oriented counter-clockwise and 
for each edge $\ca \in \cup_i\partial P_i$ 
there is a unique translation $\ct_{\ca} \neq 0 $ such that 
$\ca+\ct_{\ca}=-\cb \in \cup_i\partial P_i$. 
Identifying pairs of edges $\ca$ and $\ca+\ct_{\ca}$ gives a 
compact surface $X$ which is by construction a {\em translation surface} 
with  flat metric induced by the Euclidean metric on $\C$. 
Moreover a line field in direction $\theta \in S^1$ on $\C$  
induces a foliation $\sF_{\theta}(X)$ on $X$.  
Finally the differential $dz$ descends to $X$ (vertices of the 
polygons removed) and defines a holomorphic differential $\omega$ 
on $X$ with zeros located in the vertices of the $P_i$. 
Really important is the following group operation:
\smallskip\\   
{\bf $\mathbf{SL_2(\R)}$ action on translation surfaces.}
Take the linear operation of $\slr$ on $\R^2 \cong \C$ and  
choose $A \in \slr$. 
Then the set $\cup_i AP_i \subset \C$ 
with the identification $\cup_i \partial AP_i \ni A \cdot \ca 
\leftrightarrow A \cdot \ca + A \cdot \ct_{\ca}=-A \cdot \cb 
\in \cup_i\partial AP_i$ is a translation surface $A \cdot X$, 
a deformation of $X$. 
In this way we obtain a $\slr$-action on the set of 
translation surfaces.  
\vspace*{1mm}\\
{\bf Equivalence relation.} Two translation surfaces $X$, $Y$ 
are {\em equivalent}, 
if there exists a translation diffeomorphism 
$\phi: X \rightarrow Y$, i.e. 
$\D\phi = \id$ in polygonal coordinates above.
The {\em moduli space} of equivalence classes of translation 
surfaces can be identified with the moduli space $\Omega_1 \sM_g$ of 
genus $g$ Abelian differentials with normalized area (given by 
equation \ref{unitarea}). 
\medskip\\
{\bf Elliptic differentials in genus $1$ -- Lattice surfaces.}
Take an elliptic differential $\abx$ of genus $1$, 
i.e. a Riemann surface $X$ of genus one with a holomorphic one form $\omega$. 
For simplicity we assume $\abx$ has normalized area:
\begin{equation}\label{unitarea}
\area_{\omega}(X)=\frac{i}{2}\int_X \omega\wedge \bar{\omega}=1.
\end{equation}
The {\em absolute periods}  
\[ \per(\omega):= \left\{ \int_{\gamma}\omega: \gamma \in H_1(X;\Z)
\right\} 
\subset \R^2 \cong \C\] 
of $\abx$ define a lattice in $\R^2$. 
In particular the 
elliptic differential $(\C/\per(\omega),dz)$ has the same 
absolute period lattice as $\abx$ and in natural charts 
\[z(p)=\int^p_{p_0}\omega\]
we see that locally $dz=\omega$. This in turn implies, up to orientation  
\[ \abx \cong (\C/\per(\omega),\pm dz) \leftrightarrow \per(\omega),\] 
i.e. up to sign each elliptic differential can be identified with 
a lattice $\Lambda \in \C$.  
\medskip\\
Now represent the flat torus $\R^2/\Z^2$ by the square 
$Q$ with vertices $(0,0), (1,0)$, $(1,1), (0,1) \in \Z^2$ and 
take its image $A \cdot Q \subset \R^2$ under $A \in \slr$.  
Identifying parallel sides of the parallelogram 
$A \cdot Q$ defines a new torus $\T^2_{A}$. Moreover   
the edges of $A \cdot Q$ define a lattice $\Lambda_A =A \cdot \Z^2$, 
such that   
\[A \cdot \T^2 = \T^2_{A} = \R^2/\Lambda_A.\]
It is clear that $\Lambda_A = A \cdot \Z^2 = \Z^2$ if and only if 
$A \in \slz$. 
To discover a translation homeomorphism 
\[\T^2_{A} \rightarrow \T^2 \ \mbox{ for } A \in \slz\]
think of $\R^2$ being square-tiled 
by copies of $Q$. Now consider $A \cdot Q \subset \R^2$ 
and cut $\R^2$ along the edges of the square-tiling by $Q$. 
Then reassemble $Q$ by translating the pieces obtained from  
cutting $A \cdot Q$. Note that the $A \in \slr$ 
for which vertices of $A \cdot Q$ become 
vertices of $Q$, are exactly the $A \in \slz$.      
Thus the moduli space $\sE_1$ of (normalized) elliptic 
differentials or translation tori (with a direction) equals 
\[ \sE_1 = \slr \cdot \T^2 \cong \slr/\slz. \]
 
We see that the {\em stabilizer} 
$\{A\in \slr: A \cdot \T^2 \cong \T^2\} \cong \slz$ 
is a {\em lattice} in $\slr$. For general $\abx$  
one writes 
\[\slv{X,\omega}:=\{A\in \slr: A \cdot \abx \cong \abx\}\]
and calls this stabilizer the {\em Veech group} of $\abx$.
It is already remarkable that there are  
translation surfaces $\abx$ with a lattice 
stabilizer $\slv{X,\omega} \subset \slr$ which are not 
elliptic differentials. Translation surfaces with 
lattice stabilizer are called {\em lattice surfaces} 
or {\em Veech surfaces} in honor of 
W. Veech who found the first series of 
lattice surfaces, which are not elliptic differentials, by unfolding 
billiards in the regular $n$-gon \cite{v2,v3}. 
Recently C. T. McMullen \cite{mcm1} and Kariane Calta \cite{c} 
discovered that L-shaped polygons 
with a certain algebraic condition on the length of their sides 
are also lattice surfaces. Infinitely many Veech surfaces  
constructed from L-shaped tables, or from regular $n$-gons, 
are {\em not} elliptic differentials.
\medskip\\
{\bf Action of SL$\mathbf{(X, }\ \mathbf{\omega)}$ on coverings.} 
Given a lattice surface $\abx$ and a branched covering 
$\pi: (Y,\tau=\pi^{\ast}\omega) \rightarrow \abx$, a deformation of $\aby$ by 
$A \in \slv{X,\omega}$ is again a cover of $\abx$. 
Thus we get an operation of $\slv{X,\omega}$ on coverings 
of $\abx$   
\medskip\\
The set of all coverings of $\T^2$ of {\em fixed genus,  
fixed degree and branch points of fixed order} has  
a natural manifold structure. Moreover $2$-dimensional 
components of this {\em modular fiber} 
are elliptic differentials and cover $(\T^2,dz)$. 
We denote $2$-dimensional modular fibers by $(\sF,\omega_{\sF})$ 
or short $\sF$. 
\smallskip\\
The goal of this paper is to describe a method relating  
asymptotic constants of an elliptic differential $\abx$
to the translation geometry and topology 
of the modular fiber $\sF$ containing $\abx$. 
This approach to elliptic differentials can be easily extended 
to covers of lattice- or Veech-surfaces.    
We also give examples of modular fibers $\sF$ 
and an easy calculation of asymptotic constants 
using our formula.
The asymptotic constants of $\abx \in \sF$ depend on the 
{\em translation geometry} of $(\sF, \omega_{\sF})$ 
and on the orbit closure $\overline{\slz \cdot \abx} \subset \sF$.    
Since it is sometimes hard to characterize a connected component 
of the modular fiber by topological and geometrical invariants 
we will follow another way.
\medskip\\
{\bf The modular fiber.} 
Assume $\abx$ is a lattice surface with lattice group $\slv{X,\omega}$
We define a space of coverings $\sF_{\omega,\tau}$ using a given 
cover $\aby \rightarrow \abx$ and the $\slv{X,\omega}$ action 
on covers of $\abx$ as 
\[\sF_{\omega,\tau} :=\overline{\slv{X,\omega}\cdot \aby}.\]
There are two possible ways of taking a closure of $\slv{X,\omega}\cdot \aby$: 
\begin{itemize}
\item[1.] inside the space of differentials with fixed 
number and orders of zeros, or 
\item[2.] including all limiting surfaces, thus  
{\em degenerated} surfaces appear in $\sF_{\omega,\tau}$.  
\end{itemize} 
For this paper we will assume $\sF_{\omega,\tau}$ 
is obtained with respect to the {\em first closure}.   
In our particular cases it is not hard to see 
that $\sF_{\omega,\tau}$ is always an open complex space. 
If $\aby \rightarrow \abx $ has $n$-branch points 
one can show \cite{s2}, that there is a natural map 
\[\pi_{\ast}:\sF_{\omega,\tau} \rightarrow X^n.\] 
This map is either a covering of $X^n$, or a covering of 
an $\slv{X,\omega}$-invariant subspace of $X^n$.    
\medskip\\
To avoid technical difficulties, we only discuss  
covers of $\T^2$ branched over exactly $2$ named points. 
Tracking the two branch points on the base torus 
while deforming a cover in the modular fiber 
gives a map to $\T^2 \times \T^2-\{[x,x]: [x] \in \T^2\}$. 
Since $\R^2/\Z^2$ acts by translations 
on $\T^2$ and on the modular fiber
we divide out this torus action. Equivalently   
we might assume one of the branch points 
is fixed, say at $[0] \in \T^2$ and obtain a covering map 
$\sF_{\tau} \rightarrow \T^2-\{[0]\}$. 
Here we have simplified $\sF_{\tau}:=\sF_{\tau,\omega}$, 
because $\omega=dz$. Now the translation structure of $\T^2$ 
pulls back to $\sF_{\tau}$ and we obtain an elliptic differential    
\[(\sF_{\tau}, \omega_{\tau}):= 
(\sF_{\tau}, \pi^{\ast}dz)\] 
which by $\slz$ invariance is a 
union of lattice surfaces with Veech group $\slz$.
Now we can take the {\em unique compactification} $\sF^c_{\tau}$ 
of $\sF_{\tau}$ which makes the continuation of $\omega_{\tau}$ 
to $\sF^c_{\tau}$ {\em holomorphic}. 
\bigskip\\
{\bf Degenerated translation surfaces.}
To understand the geometry of $\sF_{\omega,\tau}$ one 
needs to look at degenerated surfaces $X_{deg}$. These are  
just deformed Abelian differentials obtained 
by moving two or more cone points into one point. There are cases when the  
degeneration process leads to a union of two or more translation 
surfaces, which are connected in some special points only. 
In algebraic geometry degenerated surfaces are known as 
{\em stable, nodal curves}. 
\medskip\\
{\bf Example:} The space $\sF_d(1,1)$ of elliptic differentials 
$\abx$ with {\em two distinguished zeros} $z_1 \neq z_2$ of {\em order} 
$1$ each, $\per(\omega)=\lat$ and $\deg(\pi)=d$, 
$\pi: \abx \rightarrow \C/\per(\omega)$ covers $\T^2-\{[0]\}$ 
and carries a natural Abelian differential $\omega_d$. 
If we take a differential in $\sF_d(1,1)$ and collapse its two cone 
points into one, we obtain a surface with 
one cone point of order $3$, or a degenerated differential. 
Differentials with one order $3$ cone point in turn are 
order $3$ cone points of $(\sF_d(1,1),\omega_d)$.  
In \cite{ems} we show that $\omega_d$ has exactly 
\begin{equation}\label{conecount}
|Z(\omega_d)|=\frac{3}{8}(d-2)d^2\prod_{p|d}(1-\frac{1}{p^2})
\end{equation}
 zeros, all of order $2$. There are other 
surfaces in the compactification of $\sF_d(1,1)$, which are 
degenerated in the sense above: these surfaces are either 
two tori identified in one point, or a torus with with two points 
identified. The total number of degenerated surfaces $X_{deg}$
in $\sF_d(1,1)$ is 
\begin{equation}\label{degcount} N_{deg}(d)
=\frac{1}{24}(5d+6)d^2\prod_{p|d}(1-\frac{1}{p^2})\quad 
\mbox{ for } d \geq 3
\end{equation}
and $N_{deg}(2)=4$.  This is the order of 
a union of $\slz$-orbits on $\sF_d(1,1)$. 
For the rest of the paper we  
use the {\em Euler $\varphi$ function}
and the {\em Dedekind $\psi$ function}:
\begin{equation}
\varphi(d):=d \prod_{p|d}\left(1-\frac{1}{p} \right), \quad 
\psi(d):= d \prod_{p|d}\left(1+\frac{1}{p} \right)
\end{equation} 
to write 
\[ d^2\prod_{p|d}(1-\frac{1}{p^2})=\varphi(d)\psi(d). \] 
{\em Remark.} The counting formul{\ae} \ref{conecount} 
and \ref{degcount} were independently 
discovered by Kani \cite{ka1,ka2} with motivation and 
tools from algebraic geometry.
\bigskip\\




\section{Results and applications} 

With the conventions of the previous example, we establish  
\begin{theorem}\label{connect}
The modular fiber $(\sF_d(1,1),\omega_d)$ is connected.  
In particular the Veech group 
of  $(\sF_d(1,1),\omega_d)$ is 
\[ \slv{\sF_d(1,1),\omega_d} \cong \slz. \] 
\end{theorem}
{\bf Remark.} Connectedness of $\sF_d(1,1)$ was already established by 
W. Fulton \cite{fu}, as the author learned from C. T. McMullen.  
However at the end of the paper 
we prove connectedness of $\sF_d(1,1)$ by  
using that it is a $\slz$-orbit of a 
torus-cover in moduli space.  
\smallskip\\
All modular fibers $\sF_d(1,1)$ admit an involution $\sigma$  
with linear part $-\id \in \slz$, thus we can consider the  
the double cover 
\[\pr_{\sigma}:\sF_d(1,1) \rightarrow \sF_d(1,1)/\sigma\] 
The quotient {\em quadratic differential}  
$(\sF_d(1,1)/\sigma, q_{d})$, or short $\sF_d(1,1)/\sigma$,  
parameterizes (normalized) degree $d$ elliptic differentials $\aby$ 
with two {\em un}-distinguishable cone-points of order $1$. 
Now any $\aby \in \sF_d(1,1)/\sigma$ admits a hyperelliptic involution, 
which interchanges its two cone points. In particular: 
distinguishing cone points destroys the hyperelliptic involution of $\aby$. 
\begin{kor}\label{genus}
We have:
\begin{equation}
\begin{split}
&\chi(\sF^c_d(1,1))=-\frac{3}{4}(d-2)\varphi(d)\psi(d)\quad 
\mbox{for }\ d \geq 2  \ \mbox{ and }\\
&\chi(\sF^c_d(1,1)/\sigma)=-\frac{1}{12}(d-6)\varphi(d)\psi(d)
\quad \mbox{ for }\ d \geq 3, \ \mbox{while }
\end{split}
\end{equation} 
$\chi(\sF^c_2(1,1)/\sigma)=2$. 
In particular $\sF^c_d(1,1)$ is hyperelliptic, if and only if 
$d=2,3,4$ and $d=5$. The surface $\sF^c_d(1,1)/\sigma$ is a torus if and 
only if $d=6$.

Moreover the parity $\Psi$ of the spin structure defined by the 
meromorphic quadratic differential $q_{d}$ on 
$\sF_d(1,1)/\sigma$ is 
\begin{multline}\label{spin}
\Psi(q_d)= \frac{|\chi(\sF^c_d(1,1)/\sigma)|}{2}\! \mod 2 
\equiv 
\\
 \equiv \left\{\begin{array}{ll} 1 &\mbox{ if } \ d=2,3,4,5 \\
0  &\mbox{ if } \ d =2n \geq 6 \\
\frac{1}{24}\varphi(d)\psi(d)\! \mod 2  &\mbox{ if } \ d =2n+1 \geq 7. 
\end{array}\right.
\end{multline} 
\end{kor} 
{\bf Remark.} The proof of these identities appears at the end of the paper. 
Let $\CH:= \{z\in \C: \im(z)>0 \}$ be the 
{\em Poincar\'e upper half plane} and 
$\Gamma(d)$ the {\em principal congruence subgroup of level} $d$.
In \cite{ka1,ka2} E. Kani describes the Hurwitz-Scheme of genus 
$2$ elliptic differentials and found (over $\C$) a description as 
an open subscheme of 
\[ X(d):= \Gamma(d)\backslash \CH \quad !\] 
Here the quotient space   
$X(d)$ is obtained by considering the action of $\slz$ 
by rational transformations on $\H$. The advantage of 
looking at the affine linear action of $\pslz$ on  
$\sF_d(1,1)/\sigma$ is, that it commutes with 
the $\pslz$-action on surfaces parameterized by 
$\sF_d(1,1)/\sigma$. 
\medskip\\
The results above allow to describe asymptotic 
quadratic growth constants in terms of the modular fiber. 
Quadratic growth rates are best expressed 
in terms of a {\em Siegel-Veech constant} \cite{v4, em, emz}: 
\[\frac{\pi}{\zeta(2)} \cdot c_{cyl}(\alpha):=
\lim_{T \rightarrow \infty}\frac{N(Cyl(\alpha),T)}{T^2}\]
where 
\[N(Cyl(\alpha),T):=|\{hol(c) \subset \R^2: c \in Cyl(\alpha)\}
\cap \{(x,y) \in \R^2: x^2+y^2 \leq T^2 \}| \] 
Here $\hol(c)=\int_{\gamma}\alpha$ and $\gamma$ is any geodesic 
loop around the core of $c$.
\begin{theorem}[{\bf Cylinders}]\cite{s2,s3} \label{siegelveech}
Assume $(\sF_{\tau}, \alpha_{\tau})$ is a  
modular fiber parameterizing elliptic differentials $\abas$ with exactly $2$ 
cone points. Suppose the     
horizontal foliation of $\sF_{\tau}$ 
decomposes into open cylinders 
$\sC_1,...,\sC_{n_{\sF}}$ of periodic regular leaves, 
bounded by singular leaves 
$\partial^{top}\sC_1,...,\partial^{top}\sC_{n_{\sF}}$. 
Then every elliptic differential 
$\abas \in \sC_i$ 
has a completely periodic horizontal foliation 
where the number of cylinders, say $n_i$, depends only on $\sC_i$. 
The horizontal cylinders of $\abas$ have width $w_{i,1},...,w_{i,n_i}$ 
independent of $\abas \in \sC_i$. If $\abas \in \sF_{\tau}$ 
has infinite $\slz$ orbit, its Siegel-Veech 
constant is: 
\begin{eqnarray}\label{genericcount}
c_{cyl}(\alpha)=
\frac{1}{\area(\sF_{\tau})}\sum^{n_{\sF}}_{i=1}
\sum^{n_{i}}_{k=1}\frac{\area(\sC_i)}{w^2_{i,k}}.
\end{eqnarray} 
If the $\slz$ orbit 
$\sO_{\alpha}:=\{A \cdot \abas:\ A \in \slz \} \subset \sF_{\tau}$ 
of $\abas$ is finite, we have 
\begin{eqnarray}\label{torsioncount}
c_{cyl}(\alpha)=
\frac{1}{|\sO_{\alpha}|}
\sum^{n_{\sF}}_{i=1}\left(  
\sum^{n_{i}}_{k=1}\frac{|\sO_{\alpha} \cap \sC_i|}{w^2_{i,k}}+  \right. 
\left. \sum^{m_i}_{k=1}
\frac{|\sO_{\alpha} \cap \partial^{top} \sC_i|}{w^2_{i,k}}
\right).
\end{eqnarray}
{\em Remark.} Connectedness of the modular fiber is not necessary 
to obtain this Theorem. Note however that $(\sF_{\tau},\alpha_{\tau})$ 
is a union of surfaces if it is not connected.
As for differentials contained in $\sC_i$, 
the number and width of cylinders contained 
in the horizontal foliation of $\abas \in \partial^{top}\sC_i$ 
depends on $\partial^{top}\sC_i$ only. 
\end{theorem}
Since asymptotic constants for cylinders on  
$\abas \in \sF_{\tau}$ depend on the cylinder-
decomposition of $\CF_h(\sF_{\tau})$, it 
is not surprising that the asymptotic constants 
for saddle connections connecting the two different 
cone points of $\abas$ depend on the {\em saddle connections} 
of $\CF_h(\sF_{\tau})$.  Note, that we consider 
{\em degenerated surfaces} in the closure of 
$\sF_{\tau}$ as {\em marked} points of $\sF_{\tau}$ 
and therefore find typically more saddle connections in 
$\CF_h(\sF_{\tau})$ as we expect from recognizing only 
cone points. 

Denote the set of saddle connections 
contained in $\CF_h(\sF_{\tau})$ by $SC_h(\sF_{\alpha})$ 
and note that {\em as a set of singular leaves} 
\[SC_h(\sF_{\alpha})  = \bigcup^n_{i=1} \partial^{top} \sC_i.\] 
Let us take $\abas \in s \in SC_h(\sF_{\alpha})$ and  
deform $\abas$ along $s$ into the right ($+$) or left ($-$) 
endpoint of $s$ in $\sF^c_{\alpha}$. That means we degenerate $\abas$ 
into a cone point or a point representing a degenerated 
surface of $\sF^c_{\alpha}$.
Tracking the family of deformed surfaces we see that 
along the deformation of $\abas$ we degenerate $m^{+}_s$ 
($m^{-}_s$) horizontal saddle connections of length 
$s^{+}_{\alpha}$ ($s^{-}_{\alpha}$ respectively) on $\abas$.  
Note that $s^{\pm}_{\alpha}$ equals the distance of 
$\abas \in s$ to the right ($+$) or left ($-$) 
endpoint of $s$. 
\medskip\\
If $o_1$ and $o_2$ are the orders of the two zeros of $\alpha $ 
then the maximal number $m$ of saddle connection which 
can be killed by one deformation is $\min(o_1,o_2)$.
To keep the following statement as elementary as possible, 
we name all the degenerated points and cone points of 
$(\sF^c_{\tau},\omega_{\tau})$ and assume 
the list is given by $z_1,...,z_{n_{\tau}}$. Associated to  
this list we get a list $o_1,...,o_{n_{\tau}}$ of {\em orders} 
of the $z_i$ and a list $m^+_1,...,m^+_{n_{\tau}}$ 
of {\em multiplicities}, telling us how many saddle connections 
disappear while degenerating a surface into $z_i$ from the 
right along a horizontal saddle connection. 
By walking along a small circle around the (cone-)point 
$z_i \in \sF^c_{\tau}$ one can 
see that $m^+_i$ is well-defined, i.e. the same for 
each horizontal saddle connection $s$ terminating in $z_i$. 
\medskip\\
Before we state the Theorem, we like to mention that 
there are asymptotic constants (see \cite{s3}) 
which reflect finer properties of $\abas$ and $\sF_{\tau}$, 
for instance one can use different weights $m^{\pm}_s$ 
associated to topological/geometrical 
properties of the surfaces represented 
by the special points $z_i \in \sF_{\tau}$ (see \cite{s3}).  
One can restrict to certain subsets of the set of 
cone points or the set of horizontal saddle connections in 
$\sF_{\tau}$ too.

\begin{theorem}[{\bf Saddle connections}]\cite{s2,s3}
With the assumptions and notations of Theorem \ref{siegelveech}, 
we find for the asymptotic quadratic growth rate $c_{\pm}(\alpha)$ 
for saddle connections on $\abas \in \sF_{\tau}$ connecting 
the two different cone points of $\abas$ 
\begin{eqnarray}\label{torsadcount}
c_{\pm}(\alpha)=
\frac{2}{|\sO_{\alpha}|}\sum_{s \in SC_h(\sF_{\alpha})}
 \sum_{(Z,\nu) \in \sO_{\alpha}(s)}\frac{m^+_s}{(s^{+}_{\alpha})^2}.
\end{eqnarray}
with $\sO_{\alpha}(s):=\sO_{\alpha} \cap s$ in the finite orbit case. 
For generic $\abas$ we find 
\begin{eqnarray}\label{gentorsadcount}
c_{\pm}(\alpha)=
\frac{2\zeta(2)}{\area(\sF_{\alpha})}\sum^{n_{\tau}}_{i=1}
m^+_i\widehat{o}_i \quad \mbox{ where } \widehat{o}_i=o_i+1. 
\end{eqnarray}
\end{theorem}
{\bf Remark.} The straightforward generalization 
of the above Theorem \cite{s3} includes:    
modular fibers of {\em higher dimension}, 
{\em arbitrary lattice group} $\slv{X,\omega}$ 
and {\em disconnected fibers} $\sF_{\tau,\omega}$. 
\smallskip\\
Depending on the specific problem, 
the formul{\ae} presented here tie the evaluation 
of Siegel-Veech constants of an elliptic differential 
$\abx \in \sF$ to
\begin{itemize} 
\item the counting of certain types 
degenerated surfaces
in the closure of the modular fiber $\sF$ 
\item the counting of cone points of $\omega_{\sF}$   
\item the classification of finite $\slz$-orbits in $\sF$.\\ 
To determine  
the constants for saddle connections on lattice elliptic differentials 
one needs to know 
\item the intersection of a particular $\slz$ orbit 
with $SC_h(\sF)  = \bigcup^n_{i=1} \partial^{top} \sC_i$.
\vspace*{3mm}
\end{itemize}

{\bf A basic example.} To apply the whole method we take the 
example of two marked tori, worked out by the author in  
\cite{s1}. Take the torus $\T^2=\C/\lat \cong \R^2/\Z^2$ 
marked in two points. We assume one of the marked points is 
$[0]:= 0 + \Z^2 $, if the other is $[m] \neq [0]$ we write 
\[\T^2_{[m]}= (\C/\lat,[0],[m]) \cong (\R^2/\Z^2,[0],[m]) \]   
The moduli space of $2$-marked tori is simply the torus 
$\T^2- \{[0]\}$ if we agree to distinguish the marked points. 
\medskip\\
Now the horizontal foliation of $\T^2- \{[0]\}$ 
consists of one cylinder $\sC$ (of height and width one) and 
one saddle connection $\partial^{top}\sC$ connecting $[0]$ 
with itself. Now the horizontal foliation of the torus $\T_{[m]}$ contains  
\begin{itemize}
\item two cylinders of width one if $[m]\in \sC$ and 
\item one cylinder of width one if $[m]\in \partial^{top}\sC$
\end{itemize}  
Formula \ref{genericcount} then implies that the asymptotic quadratic 
constant $c_{cyl}(gen)$ for isotopy classes of periodic trajectories for 
the generic two marked torus cover is $2$ (which is of course easy to see 
without a fancy formula). Now a surface or point in $\sF_1(0,0)$ is generic 
if and only if it has infinite $\slz$ orbit and these are exactly 
the irrational points in $\T^2 - \{[0]\}$, i.e. the 
set $\T^2 - \Q^2/\Z^2$.  
\medskip\\
{\bf Finite orbit case: torsion points on $\mathbf{\T^2}$.}
The set of torsion points of $\T^2$ is the kernel 
of the multiplication homomorphism  
\[\T^2[n]:= \ker(\T^2 \stackrel{n}{\rightarrow} \T^2)=
\frac{1}{n}\Z^2 /\Z^2,\quad \mbox{where } n: [z] \mapsto [nz].\]
It is not hard to see  that  
the $\slz$-orbit $\sO_n$ of $\frac{1}{n} \in \T^2$ is 
\begin{equation}\label{torbit} 
\sO_n =\left\{\left[\frac{a}{n}+i\frac{b}{n}\right] \in \T^2:\  
a,b,n \in \Z \mbox{ with } \gcd(a,b,n)=1 \right\}.
\end{equation} 
In particular 
\[|\sO_n|= 
n^2 \prod_{p|n}\left(1-\frac{1}{p^2} \right)=
\varphi(n)\psi(n).\] 
Thus $\slz$ operates transitively on the set $\T^2(n)$ of torsion points 
of {\em order} $n$. These are the torsion points in $\T^2[n]$ 
vanishing by multiplication with $n$, but do not vanish by 
multiplication with any $m|n$. We have  
\[\sO_n:=\slz \cdot [1/n]=\T^2(n). \]
To apply formula \ref{torsioncount} we need to know how 
many points of $\sO_n$ intersect with the line $\partial^{top}\sC$ 
and this is easy, in fact: 
\[ \sO_n \cap \partial^{top}\sC=\{[a/n] \in \T^2: \ a,n \in \Z,\ 
\gcd(a,n)=1\}\]
and thus $|\sO_n \cap \partial^{top}\sC|=\varphi(n)$.
Because of that $|\sO_n \cap \sC|=\psi(n)\varphi(n)-\varphi(n)=
\varphi(n)(\psi(n)-1)$. Now all horizontal cylinders on all the 
marked tori parameterized by $\sC$ have width $1$ and there are always 
two of them, while differentials on $\partial^{top}\sC$ admit 
only one horizontal cylinder (of width one of course).
\medskip\\ 
Altogether we find the asymptotic growth rate of periodic cylinders 
for any marked torus contained in $\sO_n$: 
\begin{equation}
c_{cyl}(n)=2\frac{\varphi(n)(\psi(n)-1)}{\psi(n)\varphi(n)}+
\frac{\varphi(n)}{\psi(n)\varphi(n)}=2-\frac{1}{\psi(n)}.
\end{equation}
Taking the limit for $n$ to infinity gives the generic constant
\[c_{cyl}(gen)=\lim_{n \rightarrow \infty}c_{cyl}(n)=2.\]
\medskip

{\bf Counting saddle connections.}  The only interesting 
question about the quadratic growth rate for saddle connections 
connecting the two different marked points. 
For each $0<k<n$ with $(k,n)=1$ we need to calculate the two distances 
of the point $[k/n] \in \T^2$ to $[0] \in \T^2$. 
Now using formula \ref{torsadcount}
we obtain for torsion points of order $n$: 
\begin{equation}
c_{\pm}(n)=2\frac{n^2}{\varphi(n)\psi(n)}\sum_{(k,n)=1}\frac{1}{k^2}.   
\end{equation}
In \cite{s1} we gave an explicit argument showing that for  
differentials parameterized by irrational (= generic) 
points on $\T^2$: 
\begin{equation}
c_{sc}(\pm)=\lim_{n \rightarrow \infty}c_{\pm}(n)=2\zeta(2).
\end{equation}

This example is the first of the series which we call 
$d$-{\em symmetric torus coverings}. 
To construct $d$-symmetric torus coverings 
one uses a connected sum construction for translation surfaces:
\bigskip\\ 
{\bf Connected sum construction.} Given an Abelian differential 
$\abx$ and a leaf $\sL \in \CF_{\theta}(X)$. Take 
$a \in \sL$ and define the line segment 
\[I:=[0,\epsilon]e^{i\theta}+a \subset \sL.\] 
Then for $d \geq 2$ and a cycle $\sigma \in S_d$ we 
define the Abelian differential
\[ (\#^d_{I,\sigma} X,\#^d_{I,\sigma} \omega )\]
by slicing $d$ named copies $X_1,...,X_d$ of $X$ 
along $I$ and identify opposite sides of the slits 
according to the permutation $\sigma$. 
The differential $\#^d_{I,\sigma} \omega$ on $\#^d_{I,\sigma} X$ 
is uniquely defined by the property 
\[\#^d_{I,\sigma} \omega|_{X_i}=\omega_i=\omega.\] 
Note: we can rename the $d$ copies of $X$ such that 
the cycle $\sigma$ becomes $\tau=(1,2,3,...,d)$. 
In this case we simply write: 
\[(\#^d_{I} X,\#^d_{I} \omega )
=(\#^d_{I,\tau} X,\#^d_{I,\tau} \omega ).\]  

If $\gamma_1$ and $\gamma_2$ are two chains of 
geodesic segments on an Abelian differential $\abx$ with 
\[\partial \gamma_1= \partial \gamma_2\]  
we might use cut and paste to see that 
\[ (\#^d_{\gamma_1}X,\#^d_{\gamma_1}\omega)=
(\#^d_{\gamma_2} X,\#^d_{\gamma_2}\omega),\] 
if $\gamma_1$ and $\gamma_2$ are isotopic along 
an isotopy fixing the endpoints 
$\partial \gamma_1$ and containing no cone points. 
Another way to say this is that the lifts of 
$\gamma_i$ to the universal covering $\widetilde{X}$ 
of $X$ bounds a disk $B$ containing no cone points 
(in its interior).  
\bigskip\\
{\bf $\mathbf{d}$-symmetric differentials.}
We apply this construction to a torus covering, 
by taking $d$ copies of $\T^2$, slice them along 
the projection of the line segment 
$I=I_v = [0,v]\subset \C$ ($v \in \C!$) 
to $\T^2$. Denote the resulting differential by 
\[(\#^d_{I}\T^2,\#^d_{I}dz).  \]
The underlying surface has genus $d$ and 
the translation structure has precisely two cone 
points of order $d$. 
We define $d$-symmetric torus coverings as follows 
\begin{itemize}
\item  $\tau$ has exactly {\em two} zeros of order $d-1$ 
\item  $\Z/d\Z \subset \aut \aby$
\item  $\deg(\pi)=\int_X \pi^{\ast}(dx \wedge dy) = d$\vspace*{1mm} 
\end{itemize} 
with the natural projection $\pi: \aby \rightarrow \C/\per(\tau) $. 
Note that all elliptic differentials of the shape 
$(\#^d_{I}\T^2,\#^d_{I}dz)$ are $d$-symmetric. 
\medskip\\
Denote the {\em set of isomorphy classes of $d$-symmetric   
coverings of } $\T^2=\C/\lat$ by $\sF^{sym}_d:=\sF^{sym}_d(d-1,d-1)$. 
Note that our previous examples, $2$-marked tori, are simply 
$1$-symmetric differentials, and $\sF_1(0,0)\cong \T^2 - \{0\}$. 
We show in \cite{s3}: 
\begin{theorem}\label{symstructure}
The set $\sF^{sym}_d$ has a natural structure as 
a torus (covering) $\T^2_d:=\R^2/d\Z^2$ without integer lattice points. 
The $\slz$ operation on $\sF^{sym}_d = \R^2/d\Z^2-\Z^2/d\Z^2$ 
commutes with the $\slz$ operation on surfaces parameterized by 
$\sF^{sym}_d$. 
\end{theorem}
For $d$-symmetric differentials we evaluate formula \ref{genericcount} 
to find the Siegel-Veech constants: 
\begin{theorem}\label{siegelveechdsymm}
Let $\abas \in \T^2_d$ be $d$-symmetric and 
$\abas \notin \Q^2/d\Z^2$, i.e. has infinite $\slz$-orbit in $\T^2_d$. 
Then the asymptotic quadratic growth 
rate of periodic cylinders $c_{cyl}$ on $S$ is: 
\begin{equation}
c_{cyl}(S)= c_{cyl}(d)=2\sum_{p|d} \frac{\varphi(p)}{p^3}.
\end{equation}
\end{theorem}
Note that by M\"obius inversion 
\[\varphi(d)=\frac{d^3}{2}\sum_{p|d} \mu\left(\frac{d}{p}\right)c_{cyl}(p).\]
Using formula \ref{torsioncount} we calculate the Siegel-Veech 
constants for $d$-symmetric 
differentials $\abas \in \Q^2/d\Z^2 \subset \T^2_d$,  
the {\em torsion points} in $\sF^{sym}_d$, as well. The various 
asymptotic constants 
depend very sensitive on the translation geometry of the 
surfaces. In particular some Siegel-Veech constants 
for special types of saddle connections are of interest. 
\medskip\\

\section{Modular fibers in genus $2$}

{\bf How to describe modular fibers.} Our approach works, 
if one is able to gain enough information about the translation 
geometry and topology of the space $\sF_{\omega,\tau}$, 
in particular one needs to count certain sets 
of cone points of the space $\sF_{\omega,\tau}$. 
We do not claim this is a trivial task, but one can do it to an 
extend making results access-able which are very hard to gain 
without using the geometry of $\sF_{\omega,\tau}$. 
\smallskip\\
For example: it takes a computer (program developed by G. Schmidthuesen 
\cite{gabi}) several days to calculate the index of 
the stabilizer of some differentials $\abas \in \sF_3(1,1)$ 
with small(!) $\slz$ orbit. On the other hand it takes 
not to long to make a picture of $\sF_3(1,1)$, using  
cylinders contained in  $\sF_3(1,1)$. In case of finite orbit 
surfaces in $\sF_3(1,1)$ it is possible to develop a formula 
for the order of the orbit \cite{s4}. 
\medskip\\
{\bf Absolute periods of $\sF_3(1,1)$.} 
The absolute period lattice $\per(\omega_{3})$ generated 
by the cylinders of $\sF_3(1,1)$ is 
\[\per(\omega_{3})=2\Z^2 \subset \Z^2.\] 
Different colors in Figure 1 show one possible 
tiling of $\sF_3(1,1)$ by squares of size $2$.
We claim that $\per(\omega_d)=2\Z^2$ for all $d \geq 2$. Here 
is an indirect argument: 
In \cite{ems} we found that $\sF_d(1,1)$ is tiled by 
$\frac{1}{3}(d-1)d\varphi(d)\psi(d)$ unit squares. 
Taking the quotient with respect to the involution $\sigma$ gives 
a map \[ \delta_{d}: \sF_d(1,1)/\sigma \rightarrow \proj^1 
= \T^2/(-\id)\] 
of degree 
\[\deg(\delta_{d})=\frac{1}{3}(d-1)d\varphi(d)\psi(d).\] 
branched over the image of $0=\T^2[1]$ under $\T^2 \rightarrow \proj^1$.
Kani \cite{ka3} on the other hand describes a map 
\[ \hat{\delta}_d: \sF_d(1,1)/\sigma \rightarrow \proj^1 \]
of degree 
\[\deg(\hat{\delta}_d)=\frac{1}{12}(d-1)d\varphi(d)\psi(d)=
\frac{1}{4}\deg(\delta_{d})\] 
which is branched over the images of the 2-torsion points $\T^2[2]$ 
under $\T^2 \rightarrow \proj^1$. 
\smallskip\\
Here is a picture of the translation surface $\sF_3(1,1)$ 
with some degenerated 
surfaces (vertices of the tiles of $\sF_3(1,1)$) shown below.  
\vspace*{1mm}\\
\begin{center}
\begin{picture}(3,0)%
\includegraphics{3fibre.pstex}%
\end{picture}%
\setlength{\unitlength}{2565sp}%
\begingroup\makeatletter\ifx\SetFigFont\undefined%
\gdef\SetFigFont#1#2#3#4#5{%
  \reset@font\fontsize{#1}{#2pt}%
  \fontfamily{#3}\fontseries{#4}\fontshape{#5}%
  \selectfont}%
\fi\endgroup%
\begin{picture}(7862,5385)(781,-5536)
\put(6001,-1186){\makebox(0,0)[lb]{\smash{\SetFigFont{9}{10.8}{\rmdefault}{\mddefault}{\updefault}{$r_{\pi/2}\sL$}%
}}}
\put(2926,-2311){\makebox(0,0)[lb]{\smash{\SetFigFont{9}{10.8}{\rmdefault}{\mddefault}{\itdefault}{$\sL$}%
}}}
\end{picture}
\begin{figure}[h]
\caption{The modular surface $\sF_3(1,1)$}
\end{figure}
\end{center}
Figure \ref{latticesurf} presents the surfaces on the
`integer lattice' of $\sF_3(1,1)$.  
\setcounter{figure}{1}
\begin{figure}[h]
\centering
\epsfig{file={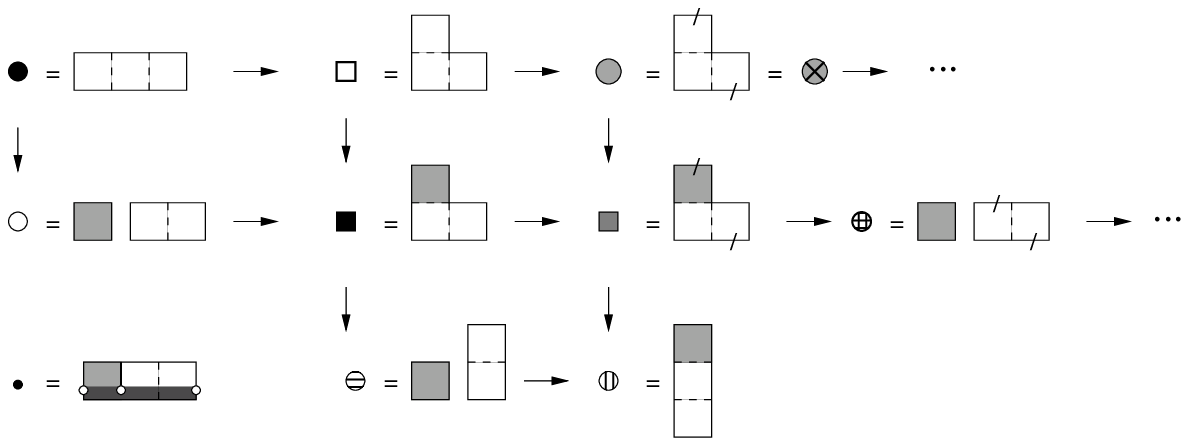},width={11cm}}
\caption{Surfaces on integer coordinates of $\sF_3(1,1)$}
\label{latticesurf}
\end{figure}
The monodromy or identification scheme of each surface in the picture 
is as follows: 
\begin{itemize}
\item {\em horizontal:} same color means same closed cylinder 
\item {\em vertical:} opposite sides are identified, unless something 
else is indicated by dashes. 
\end{itemize}
The two degenerated surfaces sitting in the middle of the slit 
in  Figure 1 are isomorphic, but appear as different points if one takes 
the closure of $(\sF_3(1,1),\omega_3)$ as elliptic differential. 

{\bf Deforming along a loop $\sL$ in $\mathbf{\sF_3(1,1)}$.}
Walking along the loop $\sL$ in $\mathbf{\sF_3(1,1)}$ from the 
left to the right represents a deformation of the degree $3$ torus 
cover denoted by the black dot to the right of the figure. 
The picture shows the $6$ surfaces at the intersection points of 
$\sL$ with the vertical edges of the tiling by squares. 
\begin{figure}[h]
\centering
\epsfig{file={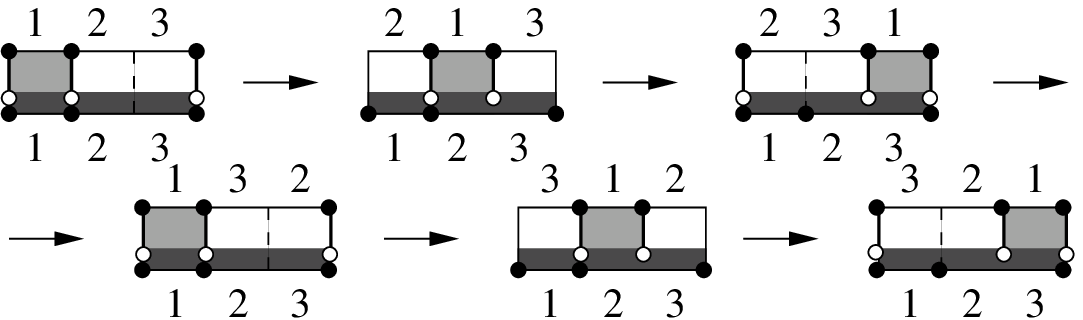},width={9cm}}
\caption{Deformation along $\sL$}
\label{deform}
\end{figure}
\vspace*{1mm}\\
Note, while deforming a surface into its neighbor, 
the vertical gluing pattern changes by a transposition. 
Note also that $\sL$ intersects $r_{\pi/2}\sL$, its image  
under rotation by $90$ degrees. 
\smallskip\\
Now we present a picture of the quadratic differential 
$q_{3}$ ($\omega^2_3=\pr^{\ast}_{\sigma}q_3$) on the sphere 
$\sF_3(1,1)/\sigma$.
\begin{figure}[h]
\centering
\epsfig{file={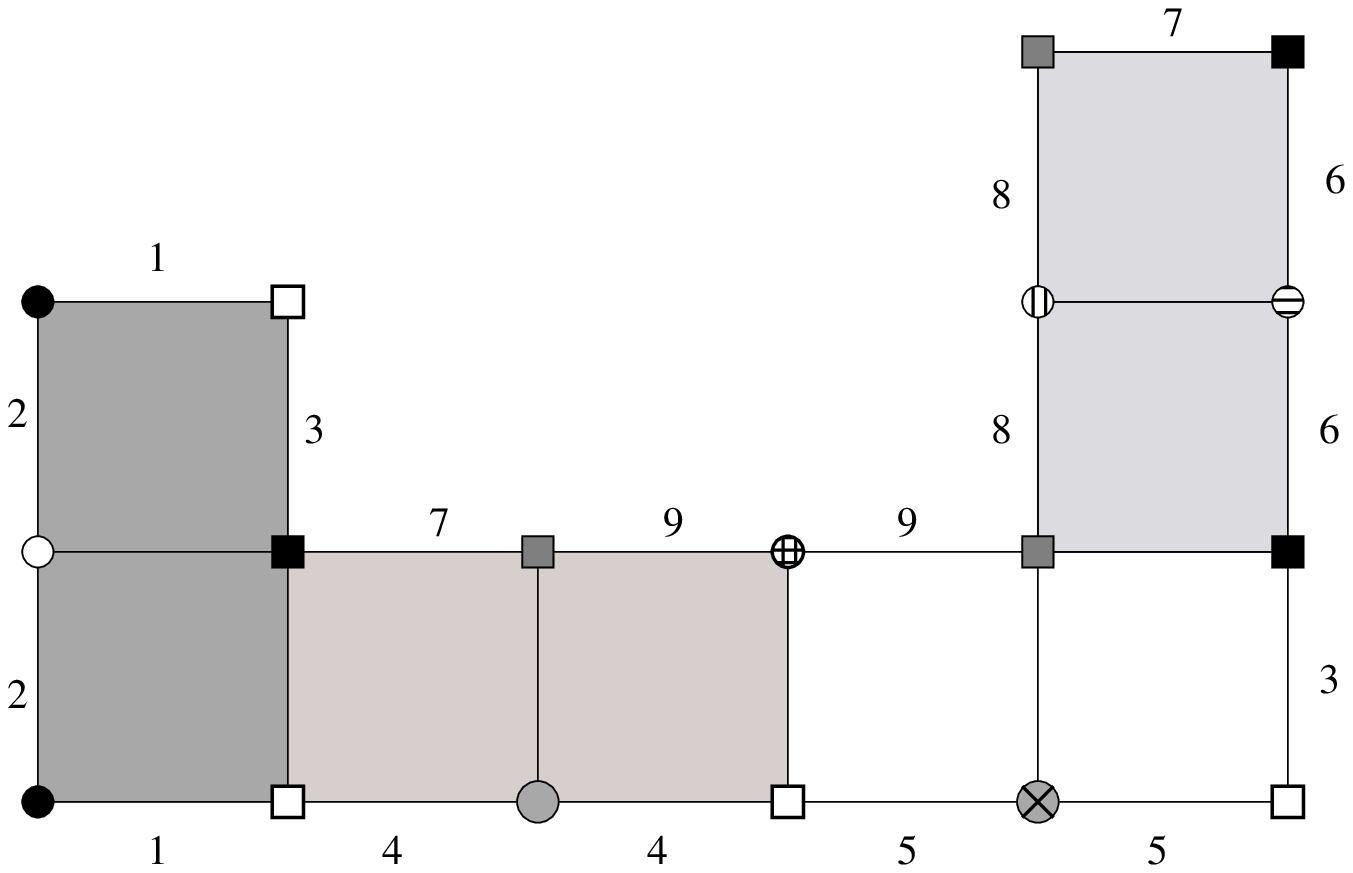},width={7cm}}
\caption{The flat sphere $\sF_3(1,1)/\sigma$}
\label{sphere}
\end{figure}
\vspace*{1mm}

{\bf Properties of $\mathbf{\sF_d(1,1)}$.} 
First we describe some translation surfaces belonging to 
$\sF_d(1,1)$. With $\T^2(a,b):=\C/a\Z\oplus ib\Z$ and a line segment 
$I=I_v:=[0,v] \subset \C$, $v \in \C$, we define the connected 
sum  
\[ \CS_{a,v}:=\T^2(a,1)\#_I\T^2(d-a,1) \in \sF_d(1,1), \mbox{ where } 
(a,d)=1.\]
The condition $(a,d)=\gcd(a,d)=1$ is necessary and sufficient 
to make sure that $\CS_{a,v}$ belongs to $\sF_d(1,1)$, and 
not to a modular fiber of lower degree $d$. 
Now assume $v=t_h+it_v$ and $t_v \in (0,1)$. 
We call $t_h$ the horizontal twist and $t_v$ the vertical twist. 
Then the horizontal cylinder decomposition of $\CS_{a,v}$ 
contains a cylinder of core width $d$ above the first cone point, 
say $z_0$, and two cylinders of width $a$ and $b=d-a$ on top 
of the second cone point $z_1=[v]$. To the  three horizontal cylinders 
we associate {\em twists}, given by $t_h(d):=t_h \mod d$ for the wide cylinder 
and by $t_h(a):=t_h \mod a$,  $t_h(b):=t_h \mod b$ respectively, 
for the narrow 
cylinders. If we pick an integer twist $t_h$, the three twists essentially 
agree with the twist part of the coordinates for $\sF_d(1,1)$, 
described in \cite{ems}.
\medskip\\
{\bf Loops and cylinders in $\mathbf{\sF_d(1,1)}$.}
By condition $(a,d)=1$ the Chinese remainder theorem implies 
the map 
\begin{eqnarray}\nonumber 
t_h & \longmapsto &(t_h(a),t_h(b),t_h(c))\\
\R &\rightarrow &\R/a\Z \oplus \R/b\Z \oplus \R/d\Z 
\end{eqnarray}
has kernel $abd\Z=a(d-a)d\Z$. From this it 
is easy to see the following 
\begin{prop}
For all $t_v \in (0,1)$ and $0<a<d$ with $(a,d)=1$, 
the map 
\[\gamma: \R/abd\Z \ni t_h \mapsto \CS_{a,t_h+it_v}\in \sF_d(1,1)\] 
is an isometrically embedded loop contained in the horizontal foliation of 
$\sF_d(1,1)$. Moreover the image of  
\[\gamma \times \id: \R/abd\Z \times (0,1) \ni (t_h,t_v) 
\mapsto \CS_{a,t_h+it_v}\in \sF_d(1,1)\]
is a maximal, horizontal cylinder $\sC^+_a \subset \sF_d(1,1)$. 
\end{prop}
{\em Remarks.} A proof of this proposition on a formal level 
requires to introduce period coordinates for $\sF_d(1,1)$ 
which we want to avoid at this place. Period coordinates are 
used in \cite{s2}.  
However it is easy to check that the loop $\gamma$ closes 
with $t_h=a(d-a)d$, if $(a,d)=1$. The cylinder $\sC^+_a$ 
is maximal because it is bounded by degenerate surfaces like 
$\CS_{a,0}$ and $\CS_{a,i}$, i.e. $t_h+it_v=0$ or $t_h+it_v=i$. 
The `$+$´ attached to $\sC^+_a$ is because of our convention 
that the cone points $z_o$ and $z_1$ are named. One obtains the cylinder 
$\sC^-_a \subset \sF_d(1,1)$ by taking 
\[\sC^-_a:=\{\CS_{a,t_h-it_v}: \ (t_h,t_v) \in \R/abd\Z \times (0,1)\}.\]

{\bf $\CU$ action on ${\mathbf \sC^+_1}$.}
To establish connectedness of $\sF_d(1,1)$ we look at the 
$\slz$ action on $\sF_d(1,1)$. In particular we 
are interested in the action of  
\[ \CU:=\left\{u_n=\left[\begin{smallmatrix} 
1 & n \\ 0 & 1 \end{smallmatrix}\right]
: n \in \Z \right\}\subset \slz \] 
on $\sC^{\pm}_1$. This looks trivial, but there is 
a nontrivial translation part caused by the $\CU$ action 
on surfaces $\CS_{1,t_h+it_v} \in \sC^{\pm}_1$. In fact we 
have 
\[u_1 \cdot \left[\begin{smallmatrix} t_h \\ t_v \end{smallmatrix}\right]=
\left[\begin{smallmatrix} t_h+t_v-d \\ t_v \end{smallmatrix}\right]. 
\] or $u_1 \cdot \CS_{1,t_h+it_v}=\CS_{1,t_h+t_v-d+it_v}$.
One can see this taking $t_v=1$, i.e. $\CS_{1,t_h+i}$, 
and using continuity of the $\slz$ action on $\sF_d(1,1)$. 
This tells us in particular that $\CU$ really {\em acts on} $\sC^{+}_1$. 
\medskip\\
Now we look to the action of the counter-clockwise rotation by $\pi/2$, 
i.e. $r_{\pi/2}:=
\left[\begin{smallmatrix} 0 & -1 \\ 1 & 0 
\end{smallmatrix}\right]\in \slz$.     
The identity 
\begin{equation} 
(r_{\pi/2} \cdot u_{d-1})\cdot \CS_{1,0}=
r_{\pi/2} \cdot \CS_{1,d-1} = \CS_{1,1-d} \in \partial \sC^{+}_1.
\end{equation}
shows that $r_{\pi/2} \cdot \sC^{+}_1$ 
intersects with $\sC^{+}_1$ in an open set, since 
$\CS_{1,1}$ is a smooth point. 
\medskip\\ 
Before we prove Theorem \ref{connect}, we add information on the global 
structure of the space of all torus-coverings or {\em elliptic 
covers} $\sE_d(1,1)$ of degree $d$, with 
two zeros of order one and absolute period lattice 
$\per(\omega)=\Lambda$ of covolume $1$. We have the following 
`fiber-bundle' structure: 
\begin{equation}
\sF_d(1,1) \longrightarrow \sE_d(1,1)\longrightarrow \slr/\slz.
\end{equation}
The base $\slr/\slz$ parameterizes lattices $\Lambda$ 
of covolume $1$ and $\sF_d(1,1)$ is the fiber over $\lat$. 
The other fibers are $\slr$-deformations of $\sF_d(1,1)$. 
We also need the following result
\begin{theorem}\cite{ems}\label{globconnect}
The space $\sE_d(1,1)$ is connected for all $d \geq 2$. 
In addition from each point in $\sF_d(1,1) \subset \sE_d(1,1)$ 
there is a path to the surface $\CS_{a,i\epsilon}$ 
for an $0<a<d$ with $(a,d)=1$. In particular $\sF_d(1,1)$ 
admits at most $\varphi(n)/2$ connected components. 
\end{theorem}
Now we can prove Theorem \ref{connect}:
\begin{proof}
Assume $\sF_d(1,1)$ is not connected. Since 
$\sE_d(1,1)$ is connected by Theorem \ref{globconnect}, 
all components of $\sF_d(1,1)$ must 
be on a single $\slz$ orbit. In particular there is an 
affine map of $\sF_d(1,1)$, induced by the $\slz$ action, 
which permutes the components of $\sF_d(1,1)$. Now  
$\slz$ is generated by $u_1$ and $r_{\pi/2}$ and  
$u_1$ fixes $\sC^+_1$, therefore it stabilizes a component of $\sF_d(1,1)$. 
Because  $r_{\pi/2} \sC^+_1\cap \sC^+_1 \neq \emptyset$, 
$r_{\pi/2}$ stabilizes the same connected component 
and the statement follows.  
\end{proof}
{\bf Remark 2.} The above is a relatively simple 
strategy to show connectedness of fibers $\sF$.   
To recall, take a loop $\sL$ in the modular fiber and prove that it is 
stabilized by the parabolic map \[u_v= \left \{g \cdot
\left[\begin{smallmatrix}1 & 1 \\ 0 & 1\end{smallmatrix}\right]
\cdot g^{-1}:\ g \in \slz \mbox{ with }\ g
\left[\begin{smallmatrix}1 \\ 0 \end{smallmatrix}\right] =v \right \} 
\in \slz\] 
fixing the holonomy-image $v=\hol(\sL) \in \R^2$. Then show that the 
fiber containing $\sL$ is stable under rotation by 
$r_{\pi/2}$. The method applies well in case the fiber 
of elliptic differentials is an orbit closure:
\[\sF_{\alpha}:= \overline{ \slz \cdot \abas}
\subset \overline{ \slr \cdot \abas} = 
\sE_{\alpha}.\]  
The argument is in general not sufficient, if the 
space of elliptic differentials is obtained 
by fixing algebraic or topological invariants of differentials. 
For example $\sE_d(1,1)$ is given by taking all differentials $\abx$ 
of degree $d$, i.e. with canonical map $ X \rightarrow \C/\per(\omega)$   
of degree $d$, and $\omega$ has precisely two zeros of 
order $1$. In this case one needs to establish 
connectedness of the whole space $\sE_d(1,1)$ first, see \cite{ems}.    
\bigskip\\
{\bf Finite $\mathbf{SL_2(\Z)}$ orbits in $\mathbf{\sF_3(1,1)}$.} 
The next and final step is to classify all finite $\slz$-orbits 
contained in $\sF_3(1,1)$. We will address this more 
generally in \cite{s4}. After this is done the 
asymptotic formul{\ae} can be evaluated if one is able to 
count how many points on each orbit are contained in 
each horizontal cylinder of $\sF_3(1,1)$. 
Since the generic constant for cylinders of periodic trajectories 
only depends on the horizontal cylinder decomposition of 
$\sF_3(1,1)$, we can easily evaluate the generic asymptotic 
constant for $\sF_3(1,1)$ and find:   
\begin{equation}c_{cyl}(gen)=\frac{1}{16}
\left[12\left(\frac{1}{1}+\frac{1}{2^2}+\frac{1}{3^2}\right)+
4\left(\frac{2}{1}+\frac{1}{2^2}\right)\right]=\frac{19}{12}.
\end{equation}
In large $d$, the counting of horizontal cylinders in $\sF_d(1,1)$ 
is non-trivial, see \cite{ems} for a coordinate approach. 
\medskip\\ 
{\bf Proof of Corollary \ref{genus}.} 
Suppose $\abx$ is an Abelian differential with (named) zeros $z_i$ of 
order $o_i$, then the {\em Gauss-Bonnet formula} 
for translation surfaces says 
\begin{equation}
\chi(X)= 2-2g(X)=-\sum^n_{i=1}o_i. 
\end{equation}
For quadratic differentials $(Y,q)$ (with simple poles)  
there is a similar formula 
\begin{equation}
2\chi(Y)= 4-4g(Y)=n_{-1}-\sum^n_{i=1}o_i, 
\end{equation}
where $n_{-1}$ is the number of poles of $q$ of order $1$ 
and $o_i$ is the order of the $i$-th zero $z_i$ of $q$. 
A zero (pole) of order $o_i$ is a cone point of total angle 
$(o_i+2)\pi$ w.r.t. the {\em half-translation} structure on $(Y,q)$.
\smallskip\\
The expression for $\chi(\sF^c_d(1,1))$ 
($d \geq 3$) comes from the 
fact that $\sF^c_d(1,1)$ has $\frac{3}{8}(d-2)\varphi(d)\psi(d)$ 
cone-points, all of order $3$. These cone-points are order $2$ 
zeros of $\omega_d$.  

To calculate $\chi(\sF^c_d(1,1)/\sigma)$ we note that $q_d$ 
has $n_{+1}=\frac{3}{8}(d-2)\varphi(d)\psi(d)$ cone-points with 
total angle $3\pi$, these are simple zeros of $q_d$. 
The number $n_{-1}$ of simple poles of $q_d$ equals the number of 
cone points of total angle $\pi$ on $\sF^c_d(1,1)/\sigma$, 
which in turn equals the number $N_{deg}(d)$ of {\em degenerated 
surfaces} in $\sF^c_d(1,1)$. Thus we find the stated expression from 
 \begin{equation}
\begin{split}
&2\chi(\sF^c_d(1,1)/\sigma)= N_{deg}(d)-|Z(q_d)|= n_{-1}-n_{+1}=\\
&=\frac{1}{24}\left((5d+6)-9(d-2)\right)\varphi(d)\psi(d)=
-\frac{1}{6}(d-6)\varphi(d)\psi(d)
\quad \mbox{ for }\ d \geq 3.
\end{split}
\end{equation}     
For $d=2$ we have $\chi(\sF^c_2(1,1)/\sigma)=\chi(\T^2/(-\id))=
\chi(\proj^1)=0$. 
\medskip\\
Since $q_d$ has only simple poles and zeros of order $1$, 
Theorem 1.2 in \cite{l} gives formula 
\ref{spin}, after observing that   
\[ \frac{n_{-1}-n_{+1}}{4}= 
\frac{\chi(\sF^c_d(1,1)/\sigma)}{2}= -\frac{1}{24}(d-6)\varphi(d)\psi(d) 
\in \Z \quad \mbox{ for } d \geq 3\] 
and obvious simplifications when taking this expression modulo $2$.
\qed
\bigskip\\
{\bf Acknowledgments.} The author is grateful for the   
atmosphere during the Activity on Algebraic and Topological Dynamics 
in Bonn. In particular I have to thank Sergij Kolyada  
for the excellent organization and Don Zagier for an additional short 
term invitation.   
While part of the activity I enjoyed discussions about the subject 
of this note with Anton Zorich, Curt McMullen, Richard Schwartz 
and Pascal Hubert, just to name a few.     
The author thanks the referee for many useful suggestions and 
for a careful reading of this partly overview article. 
Since some of the cited results, particularly the results about $d$-symmetric 
differentials, were not peer-reviewed while this paper was under 
review, all mistakes are entirely the fault of the author. 



\begin{thebibliography}{99}
\bibitem[C]{c}Calta, Kariane Veech surfaces and complete periodicity 
in genus 2. J. Amer. Math. Soc. 17 (2004), no. 4, 871--908 (electronic).
 \bibitem[Dij]{dijk}Dijkgraaf Robbert. Mirror symmetry and elliptic curves. The moduli space of curves (Texel Island, 1994), 149--163, Progr. Math., 129, Birkhäuser Boston, Boston, MA, 1995. 
\bibitem[EMcM]{emcm} Eskin, Alex; McMullen, Curt Mixing, counting, and 
equidistribution in Lie groups.  Duke Math. J.  71  (1993),  no. 1, 181--209.
\bibitem[E98]{e98} Alex Eskin,  Counting problems and semisimple groups,
 Doc. Math. J. DMV , Extra Volume ICM II (1998) 539-552.
\bibitem[EM98]{em} A. Eskin, H. Masur; Pointwise asymptotic
formulas on flat surfaces,  Ergodic Theory Dynam. Systems  
21  (2001),  no. 2, 443--478.
\bibitem[EMM]{emm} A. Eskin, J. Marklof, D. Morris. 
Unipotent flows and Veech surfaces. Preprint 2004.
\bibitem[EMS]{ems}Eskin, Alex; Masur, Howard; Schmoll, Martin. 
Billiards in rectangles with barriers. 
Duke Math. J. 118 (2003), no. 3, 427--463. 
\bibitem[EMZ]{emz}Eskin, Alex; Masur, Howard; Zorich, Anton.   
Moduli Spaces of Abelian Differentials: The Principal Boundary, 
Counting Problems and the Siegel--Veech Constants. 
Publ. Math. Inst. Hautes \'Etudes Sci.  No. 97 (2003), 61--179. 
\bibitem[EO]{eo}Eskin, Alex; Okounkov, Andrei. 
Asymptotics of numbers of branched coverings of a torus 
and volumes of moduli spaces of holomorphic differentials. 
Invent. Math. 145 (2001), no. 1,
  59--103.
\bibitem[EOP]{eop} Alex Eskin, Andrei Okounkov, Rahul Pandharipande. 
 The theta characteristic of a branched covering. Preprint. math.AG/0312186
\bibitem[Fu]{fu} W. Fulton. Hurwitz Schemes and moduli of curves. 
Annals of Mathematics, 90:542-575, 1969.
\bibitem[GJ]{gj}Gutkin, Eugene; Judge, Chris. 
Affine mappings of translation surfaces: geometry and arithmetic. 
Duke Math. J. 103 (2000), no. 2, 191--213.
\bibitem[GHS]{ghs}Gutkin, Eugene; Hubert, Pascal; Schmidt, Tom. 
Affine diffeomorphisms of translation 
surfaces: periodic points, Fuchsian groups and arithmeticity.  
Ann. Sci. \'Ecole Norm. Sup. (4) 36 (2003), no. 6, 847--866 (2004)
\bibitem[HL]{hl}Hubert, Pascal; Lelievre, Samuel. 
Prime arithmetic Teichm\"uller discs in $\CH(2)$. 
arXiv:math.GT\slash 0401056v2. 
To appear in Israel Journal of Mathematics. 
\bibitem[Ka3]{ka3} Kani, Ernst. Hurwitz spaces of covers of an elliptic curve. 
Survey article. http://www.mast.queensu.ca/~kani/hurwitz.htm 
\bibitem[Ka2]{ka2}  Kani, Ernst. Hurwitz spaces of genus 2 
covers of an elliptic curve. Collect. Math. 54 (2003), no. 1, 1--51.
\bibitem[Ka1]{ka1} Kani, Ernst. The number of curves of genus two 
with elliptic differentials. J. Reine Angew. Math. 485 (1997), 93--121.
\bibitem[KZo]{kzo} Kontsevich, Maxim; Zorich, Anton. 
Connected components of the moduli space of abelian differentials 
with prescribed singularities. Invent. Math. 153(2003), no. 3, 631--678.
\bibitem[L]{l} Lanneau, Erwan. Parity of the spin structure defined by a 
quadratic differential. Geometry \& Topology 8 (2004), 511-538. 
\bibitem[M4]{m4} Masur, Howard. The growth rate of trajectories 
of a quadratic differential. 
Ergodic Theory Dynam. Systems 10 (1990), no. 1, 151--176.
\bibitem[M3]{m3}Masur, Howard. Lower bounds for 
the number of saddle connections and closed trajectories of 
a quadratic differential. Holomorphic functions and moduli, 
Vol. I (Berkeley, CA, 1986), 
215--228, Math. Sci. Res. Inst. Publ., 10, Springer, New York, 1988. 
\bibitem[M2]{m2}Masur, Howard. Closed trajectories 
for quadratic differentials with an application to billiards. 
Duke Math. J. 53 (1986), no. 2, 307--314. 
\bibitem[M1]{m1}Masur, Howard. 
Interval exchange transformations and measured foliations. 
Ann. of Math. (2) 115 (1982), no. 1, 169--200. 
\bibitem[McM1]{mcm1} McMullen, Curtis T. Teichm\"uller curves 
on Hilbert modular surfaces. 
J. Amer. Math. Soc. 16 (2003), no. 4, 857-885 (electronic).
\bibitem[McM2]{mcm2} McMullen, Curtis T. Teichm\"uller geodesics of 
infinite complexity. Acta Math. 191 (2003), no. 2, 191--223.
\bibitem[McM3]{mcm3}McMullen, Curtis T. Dynamics of $\slr$ over moduli space 
in genus two, Preprint 2003.
\bibitem[McM4]{mcm4}McMullen, Curtis T. Teichm\"uller curves in genus two:
Discriminant and spin.  Math. Ann.  333  (2005),  no. 1, 87--130. 
\bibitem[McM5]{mcm5}McMullen, Curtis T. Teichm\"uller curves in genus two: 
The decagon and beyond.  J. Reine Angew. Math.  582  (2005), 173--199. 
\bibitem[McM6]{mcm6}McMullen, Curtis T. Teichm\"uller curves in genus two: 
Torsion divisors and ratios of sines. Preprint 2004.
\bibitem[GS]{gabi} Gabriela Schmidth\"usen; An algorithm for finding the 
Veech group of an origami.  arXiv:math.AG/0401185.
\bibitem[S1]{s1} M. Schmoll; On the growth rate of saddle connections 
and closed geodesics on spaces of marked tori, GAFA, 
Geom. funct. anal., Vol. 12 (2002) 622-649.
\bibitem[S2]{s2} M. Schmoll; Moduli spaces of covers of Veech surfaces II:
 Geometry and Siegel-Veech formula. In preparation. 
\bibitem[S3]{s3} M. Schmoll. Moduli spaces of covers of Veech surfaces I:
 $d$-symmetric differentials. Preprint 2004. 
\bibitem[S4]{s4} M. Schmoll. In preparation 
\bibitem[V1]{v1}Veech, William A. The Teichm\"uller geodesic flow. 
Ann. of Math. (2) 124 (1986), no. 3, 441--530. 
\bibitem[V2]{v2}Veech, W. A. Teichm\"uller curves in moduli space, 
Eisenstein series and an application to 
triangular billiards. Invent. Math. 97 (1989), no. 3, 553--583.
\bibitem[V3]{v3}Veech, William A. The billiard in a regular polygon. 
Geom. Funct. Anal. 2 (1992), no. 3, 341--379.
\bibitem[V4]{v4}Veech, William A. Siegel measures. 
Ann. of Math. (2) 148 (1998), no. 3, 895--944.
\bibitem[Vrb]{vrb} Y. Vorobets; Planar structures and billiards
  in rational polygons: The Veech alternative, Russ. Math. Surveys 51
(1996), 779-817.
\end{thebibliography}
\end{document}